\documentclass[journal,twocolumn]{IEEEtran}
\IEEEoverridecommandlockouts
\usepackage{cite}
\usepackage{amsmath,amssymb,amsfonts, mathtools}
\newcounter{relctr} %% <- counter for relations
\everydisplay\expandafter{\the\everydisplay\setcounter{relctr}{0}} %% <- reset every eq
\newcommand\labelrel[2]{%
  \begingroup
    \refstepcounter{relctr}%
    \stackrel{\textnormal{(\alph{relctr})}}{\mathstrut{#1}}%
    \originallabel{#2}%
  \endgroup
}
\AtBeginDocument{\let\originallabel\label}
\usepackage{suffix}
\usepackage{algorithmic}
\usepackage{graphicx}
\usepackage{subcaption}
\usepackage{textcomp}
\usepackage{xfrac}
\usepackage{xcolor}
\usepackage{epstopdf}
\usepackage{amsthm} 
\usepackage[justification=centering]{caption} %check if needed in two columns

\DeclarePairedDelimiterX\MeijerM[3]{\lparen}{\rparen}%
{#3\delimsize\vert\,\begin{smallmatrix}#1 \\ #2\end{smallmatrix}}

\newcommand\MeijerG[8][]{%
  G^{\,#2,#3}_{#4,#5}\MeijerM[#1]{#6}{#7}{#8}}

\WithSuffix\newcommand\MeijerG*[7]{%
  G^{\,#1,#2}_{#3,#4}\MeijerM*{#5}{#6}{#7}}

\def\BibTeX{{\rm B\kern-.05em{\sc i\kern-.025em b}\kern-.08em
    T\kern-.1667em\lower.7ex\hbox{E}\kern-.125emX}}
    
\graphicspath{{.}{../eps/}{./images/}{./images/result/}}

\newtheorem{prop}{Proposition}

\newtheorem{rem}{Remark}
\newtheorem{cor}{Corollary}
\theoremstyle{definition}

\begin{document}

\title{Statistical Properties of Transmissions Subject to Rayleigh Fading and Ornstein-Uhlenbeck Mobility}

\author{\IEEEauthorblockN{Arta Cika, Mihai-Alin Badiu, and Justin P. Coon}\\
Department of Engineering Science\\
University of Oxford, Parks Road, Oxford, UK, OX1 3PJ\\
Email: \{arta.cika, mihai.badiu, and justin.coon\}@eng.ox.ac.uk
\thanks{The authors gratefully acknowledge the support of Moogsoft Ltd. This material is also based upon work supported in part by the U. S. Army Research Laboratory and the U. S. Army Research Office under contract/grant number W911NF-19-1-0048.}}

\maketitle

\begin{abstract}
In this paper, we derive closed-form expressions for significant statistical properties of the link signal-to-noise ratio (SNR) and the separation distance in mobile ad hoc networks subject to Ornstein-Uhlenbeck (OU) mobility and Rayleigh fading. In these systems, the SNR is a critical parameter as it directly influences link performance. In the absence of signal fading, the distribution of the link SNR depends exclusively on the squared distance between nodes, which is governed by the mobility model. In our analysis, nodes move randomly according to an Ornstein-Uhlenbeck process, using one tuning parameter to control the temporal dependency in the mobility pattern. We derive a complete statistical description of the squared distance and show that it forms a stationary Markov process. Then, we compute closed-form expressions for the probability density function (pdf), the cumulative distribution function (cdf), the bivariate pdf, and the bivariate cdf of the link SNR. Next, we introduce small-scale fading, modeled by a Rayleigh random variable, and evaluate the pdf of the link SNR for rational path loss exponents. The validity of our theoretical analysis is verified by extensive simulation studies. The results presented in this work can be used to quantify link uncertainty and evaluate stability in mobile ad hoc wireless systems.
\end{abstract}

\begin{IEEEkeywords}
signal-to-noise ratio statistics, mobile ad hoc wireless networks, mobility modeling, Ornstein-Uhlenbeck process, Rayleigh fading. 
\end{IEEEkeywords}

\section{Introduction}
\IEEEPARstart{M}{obile} ad hoc networks (MANET) consist of autonomous mobile wireless devices (nodes) that can create a network in a decentralized manner, without the need for a fixed infrastructure~\cite{roy2010handbook, haas1999guest}. A link between two nodes exists if the received SNR is higher than a system-dependent threshold. In this environment, SNR changes over time due to node movements and due to variations in the propagation channel. Therefore, connections between nodes are established and broken intermittently, leading to dynamically changing network topology. Node mobility and channel randomness are thus the main factors impacting the distribution of the link SNR, and consequently, the performance of MANETs~\cite{chlamtac2003mobile, miorandi2008impact}. Particularly, routing in these systems faces strong challenges~\cite{abolhasan2004review, hoebeke2004overview}. For instance, the protocol must be able to update its routes rapidly to maintain connectivity. At the same time, it should not overflow the network with service messages. Hence, in this paper, we are motivated to investigate a fundamental issue related to node connectivity, namely, what is the probability that two mobile wireless devices are connected at any time instance? To this end, we require closed-form expressions for the distribution of the link SNR. 

The link performance is fundamentally determined by the instantaneous SNR and, hence, it has been extensively analyzed in the literature. Authors in~\cite{smith2018connectivity} used the link SNR to evaluate the mean proportion of time that two mobile devices are connected for five different types of mobility (including Ornstein-Uhlenbeck mobility). In~\cite{madadi2017shared}, the authors characterized variations in the SNR process to study the coverage and outage durations experienced by mobile users, while, in~\cite{madadi2016temporal} they studied the time variations of the SNR in the absence of fading experienced when a mobile user moves across a Poisson
cellular network. The mobility pattern used in these two papers is very simple, i.e., nodes move along straight lines in random directions with a constant speed. This model has been used extensively in the literature, but it does not fully capture the temporal dependence of the movement of a node over time~\cite{bai2003important}.

In~\cite{farkas2008link} the link SNR was proposed as a quality measure to predict links' state evolution in mesh networks. Closed-form expressions for the outage probability and the average probability of error based on the SNR of dual-hop wireless systems subject to Rayleigh fading were presented in~\cite{hasna2004performance}. The statistical properties of amplify and forward relay fading channels were studied in~\cite{patel2006statistical}. In that paper, the link SNR was used to derive the frequency of outage and the average outage duration. In the case of multi-hop scenarios, it is more difficult to find expressions for the pdf and cdf of the end-to-end SNR. In~\cite{karagiannidis2004statistical}, the authors overcame this problem by evaluating parameters such as the mean, the variance, the skewness, and the kurtosis to characterize the behavior of the distribution of the output SNR over Nakagami-m fading channels.

In this paper, we derive closed-form expressions for important statistical parameters of the link SNR and the separation distance in systems subject to both signal fading and random mobility. The SNR depends on the distance between nodes, which, in turn, is governed by the mobility model and the environmental factors controlling the channel between devices. Here, nodes move according to an OU process, a continuous-time Gaussian Markov process~\cite{gardiner2009stochastic}. This mean-reverted process is particularly suitable for describing the movement of a group of elements having the same destination suffering from random displacement around the projected trajectory. Moreover, the OU process represents a wide range of patterns with varying degrees of memory, including, as the two extreme cases, the random walk and the constant mobility model~\cite{camp2002survey}. We develop our analysis both in the absence and presence of channel randomness caused by fading. In modeling the mobile radio channel, we refer to Clarke's model~\cite{clarke1968statistical}, which assumes that the wireless channel is time-invariant during the symbol interval, but it varies spatially due to scattering.

Link connectivity is strongly dependent on the instantaneous SNR at any time $t$. Therefore, the analytical results obtained in this work can be used to create a framework to calculate the connectivity probability in closed form, and to evaluate link stability for both fading and non-fading scenarios. The main contributions of this paper are the following:
\begin{itemize}
\item We derive a complete statistical description of the squared distance process, including its pdf and autocorrelation function, and show that it forms a stationary Markov process.
\item We compute the stochastic differential equation of the link SNR. 
\item In the absence of fading, we derive closed-form expressions for the pdf, the cdf, the bivariate pdf, and the bivariate cdf of the link SNR.
\item In the presence of fading, we evaluate the pdf of the link SNR for rational path loss exponents.
\item We calculate the connectivity probability in closed form for both fading and non-fading scenarios. 
\end{itemize}
The rest of the paper is organized as follows. Section~\ref{sec:sys-model} provides basic definitions, the mobility model formulation, and the Rayleigh fading channel characterization. In section~\ref{sec:no-fading}, we present the mathematical and simulations results in mobile networks when there is no signal fading. We generalize our theoretical framework in section~\ref{sec:with-fading}, where we analyze the statistical properties of the link SNR in the presence of signal fading. Finally, concluding remarks are discussed in section~\ref{sec:conclusions}. 
%%%%%%%%%%%%%%%%%%%%%%%%%%%%%%%%%%%%%%%%
%%%%%%%%%%%%%%%%%%%%%%%%%%%%%%%%%%%%%%%%
\section{System Model}\label{sec:sys-model}
Consider two arbitrary nodes (mobile wireless devices) moving randomly over a two-dimensional plane. Each device movement is assumed to be independent from the other. The locations of the nodes at time $t \geq 0$ are given by $\mathbf{A}'_{t}=\left(X'_{t};Y'_{t}\right)$ and $\mathbf{A}''_{t}=\left(X''_t,Y''_t\right)$, respectively. We denote by $R_t$ the Euclidean distance between two nodes, $R_t = \|\mathbf{A}''_t-\mathbf{A}'_t\|$, and by $Z_t$ the squared distance, $Z_t=R_t^2$.
%%%%%%%%%%%%%%%%%%%%%%%%%%%%%%%%%%%%%%%%
%%%%%%%%%%%%%%%%%%%%%%%%%%%%%%%%%%%%%%%%
\subsection{Signal to Noise Ratio}\label{subsec:snr}
Instantaneous SNR is a critical measure that gives clear indications on the quality of a node's connectivity with its neighbors. In a typical multi-path propagation environment, the SNR of a communication link with additive Gaussian noise at any time $t$ is given by
\begin{equation}\label{eq:inst_snr}
 N_{t}=\psi Z_{t}^{-\eta/2}G_t,
\end{equation}
where $\eta$ is the path loss exponent (typically $2\leq\eta\leq5$), $\psi$ is a constant depending on different parameters such as transmit power, antenna properties, and wavelength, and $G_t$ is the channel gain. Unless specified otherwise, all results in this paper assume $\psi=1$.

A transmission from node $1$ to node $2$ at any time $t$ is successful if the SNR is higher than a certain threshold $\rho_{th}$, determined by the communication hardware, as well as the modulation and coding scheme of the wireless system. With these definitions, the probability the two nodes establish a communication link equals to
\begin{equation}\label{eq:conn_prob}
\mathbb{P}\left(Z_{t}^{-\eta/2}G_t\geq \rho_{th}\right).
\end{equation}
We are interested in quantifying the connectivity  probability~\eqref{eq:conn_prob}. To this end, we require expressions for the distribution of the link SNR.
\subsection{Mobility Model}\label{subsec:mobility}
We model the node displacements along the $x$ and $y$ coordinates, $\left\{\left(X'_t, Y'_t\right), \left(X''_t,Y''_t\right) \right\}$, by independent identical OU processes. Specifically, each coordinate $X'_t, Y'_t, X''_t,Y''_t$ is equally in distribution with $S_t$, where $\{S_t, t\geq0\}$ is an OU process. 

An OU process is defined as the solution to the stochastic differential equation (SDE)~\cite{uhlenbeck1930theory, doob1942brownian}
\begin{equation}\label{eq:sde}
\mathrm{d}S_t=\frac{1}{\tau}(\mu-S_t)\mathrm{d}t+\sqrt{D}\mathrm{d}W_t,\quad S_0=s_0\in\mathbb{R},
\end{equation}
where $\mu$ is the desired position, and $\{W_t, t\geq 0\}$ is the standard Brownian motion. The solution to the SDE in~\eqref{eq:sde} exists, and is unique~\cite{oksendal2003stochastic}. The parameters $\tau$ and $D$ are positive constants called the \textit{relaxation time} and the \textit{diffusion coefficient}, respectively; $\sqrt{D}$ controls the fluctuation in the position of the devices along each coordinate axis, and $1/\tau$ controls the rate of reversion of the device to the desired position $\mu$. Given the starting point $\{S_0=s_0, t=0\}$, the expectation is $m= \mu + (s_0-\mu)e^{-t/\tau}$ and the variance is $\alpha= \frac{D\tau}{2}\left(1-e^{-2t/\tau}\right)$. Note that $\alpha\rightarrow\frac{D\tau}{2}$ and $m\rightarrow\mu$ as $t\rightarrow\infty$.

\begin{figure}[t]
\centering
\includegraphics[width = \columnwidth]{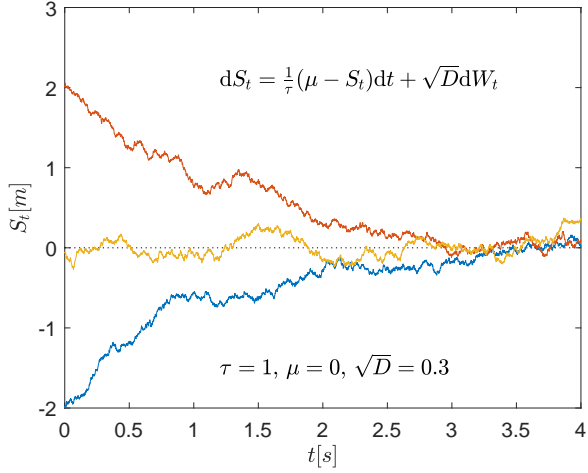}
\caption{Three sample paths of OU processes with the same parameters: $\tau = 1$s, $\mu=0$ and $\sqrt{D}=0.3$m$^2/$s, but different initial positions $S_0$.}
\label{fig:ou-process}
\end{figure}

\begin{figure}[t]
\centering
\includegraphics[width = \columnwidth]{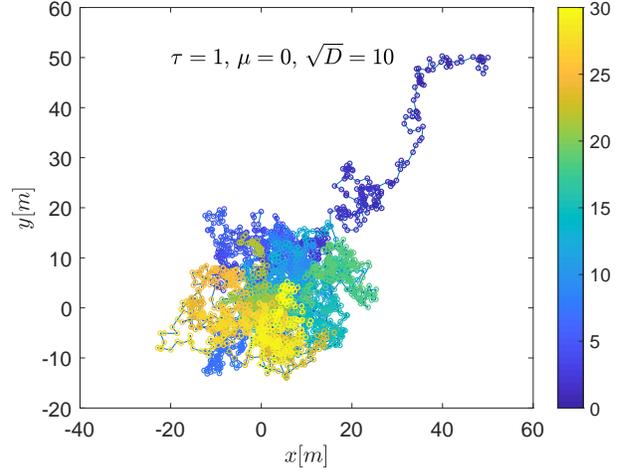}
\caption{A $2$D simulation of node trajectories along the $x$ and $y$ coordinates.}
\label{fig:ou-path}
\end{figure}

In Fig.~\ref{fig:ou-process} are shown three sample paths of OU processes with the same parameters, $\tau = 1$, $\mu=0$ and $\sqrt{D}=0.3$, but different initial positions. In the long term, all paths approach the steady-state. A $2$D simulation of node $1$ and node $2$ trajectories is plotted in Fig.~\ref{fig:ou-path}. The OU process is a continuous time Gaussian Markov process~\cite{gardiner2009stochastic}. If the initial condition of the process, $S_0$, is drawn according to the steady-state distribution, then the process is stationary. 

Another quantity of interest is the stationary autocovariance function of the OU process, which is obtained by allowing the system to approach its steady-state. It is given by~\cite{gardiner2009stochastic}
\begin{equation}\label{eq:autocorr}
K_{S}\left(\Delta t\right)=\mathbb{E}\left\{\left[S_{t+\Delta t}-m\right]\left[S_t-m\right]\right\} = \frac{D\tau}{2}e^{-\Delta t/\tau}.
\end{equation}
In the steady-state, the random variables $S_t$ and $S_u$ are only significantly correlated if $|t-u| \approx \tau$, also known as the \textit{correlation time}. 
%%%%%%%%%%%%%%%%%%%%%%%%%%%%%%%%%%%%%%%%
%%%%%%%%%%%%%%%%%%%%%%%%%%%%%%%%%%%%%%%%
%%%%%%%%%%%%%%%%%%%%%%%%%%%%%%%%%%%%%%%%
%%%%%%%%%%%%%%%%%%%%%%%%%%%%%%%%%%%%%%%%
Under the model described above, the random variables $X'_t, Y'_t, X''_t,Y''_t \sim \mathcal{N}\left(m, \alpha\right)$ are independent, and the squared distance between two nodes $1$ and $2$, at any time $t$, is given by
\begin{equation}\label{eq:distance}
Z_t= X_t^2+Y_t^2,
\end{equation}
where $X_t=X''_t-X'_t\sim \mathcal{N}\left(0, 2\alpha\right)$ and $Y_t=Y''_t-Y'_t\sim \mathcal{N}\left(0, 2\alpha\right)$ are independent random variables. 
%%%%%%%%%%%%%%%%%%%%%%%%%%%%%%%%%%%%%%%%
%%%%%%%%%%%%%%%%%%%%%%%%%%%%%%%%%%%%%%%%
\subsection{Rayleigh Fading Channel Model}\label{subsec:fading}
In modeling the mobile radio channel, we refer to Jakes' model~\cite{jakes1994microwave}. This radio model takes into account the dynamics of signal power variations that are unavoidably caused by obstructions and irregularities in the propagation path between the receiver and transmitter. The main assumption of the model is that all reflected signals at the receiver become uncorrelated in amplitude while uniformly distributed between $0$ and $2\pi$ in phase. Consequently, the in-phase and quadrature components of the channel response tend towards a Gaussian distribution, while their envelope follows a Rayleigh distribution with parameter $\lambda$. In a noisy fading channel, the most important parameter is the channel gain, $\{G_{t}, t\geq0\}$, as it determines the received power at any time instance $t$. By a simple transformation of random variables it is easy to show $G_t\sim\mathsf{Exp}\left(\lambda^2\right)$ for all $t$. In this paper, we assume that the channel gain is exponentially distributed with mean one, i.e., $\lambda=1$. It will become apparent that the omission of this detail does not hinder the development of important results.
%%%%%%%%%%%%%%%%%%%%%%%%%%%%%%%%%%%%%%%%
%%%%%%%%%%%%%%%%%%%%%%%%%%%%%%%%%%%%%%%%
%%%%%%%%%%%%%%%%%%%%%%%%%%%%%%%%%%%%%%%%
%%%%%%%%%%%%%%%%%%%%%%%%%%%%%%%%%%%%%%%%
The autocovariance function best describes the temporal variability of the channel gain. In Jakes' channel model, it is shown to be~\cite{aulin1979modified}
\begin{align}\label{eq:correl_fad}
\mathcal{K}_G\left(\Delta t\right)=\mathrm{J}_0^2\left(2\pi\nu_{\max}\Delta t\right),
\end{align}
where $\mathrm{J}_0\left(\cdot\right)$ is the zero-order Bessel function of the first kind, and $\nu_{\max}$ is the maximum Doppler shift, a measure for the rate of change of the fading channel. It follows that the channel gain $\{G_{t}, t\geq0\}$ is a wide-sense stationary (WSS) process since its mean and variance are time-invariant, and its autocovariance function depends only on the time shift $\Delta t$.
%%%%%%%%%%%%%%%%%%%%%%%%%%%%%%%%%%%%%%%%
%%%%%%%%%%%%%%%%%%%%%%%%%%%%%%%%%%%%%%%%
\section{Statistical Properties of Mobile Networks without Fading}\label{sec:no-fading}
In this section, we derive the statistical properties of the link SNR when there is no signal fading affecting the link between nodes, i.e., $G_{t}=1$ for $t\geq 0$. Neglecting the effects of the underlying wireless channel, the SNR simplifies to
\begin{align}\label{eq:snr_dist}
N_{t}=Z_{t}^{-\eta/2}.
\end{align}
Consequently, the connectivity probability can be written as $\mathbb{P}\left(Z_{t}\leq r^2_0\right)$, where $r_0 = (1/\rho_{th})^{\frac{1}{\eta}}$ defines the typical connection range and depends on several system parameters, such as the transmit power, wavelength, bandwidth, and the noise power spectral density. This model is also known as the hard connection model of link connectivity~\cite{penrose2003random}. It states that nodes can communicate whenever they lie within some critical distance of each other.  Hence, in this scenario, the statistical properties of the link SNR are completely determined by the squared distance between the two nodes. 
%%%%%%%%%%%%%%%%%%%%%%%%%%%%%%%%%%%%%%%%
%%%%%%%%%%%%%%%%%%%%%%%%%%%%%%%%%%%%%%%%
\subsection{Squared Distance Process}
\begin{figure}[t]
\centering
\includegraphics[width = \columnwidth]{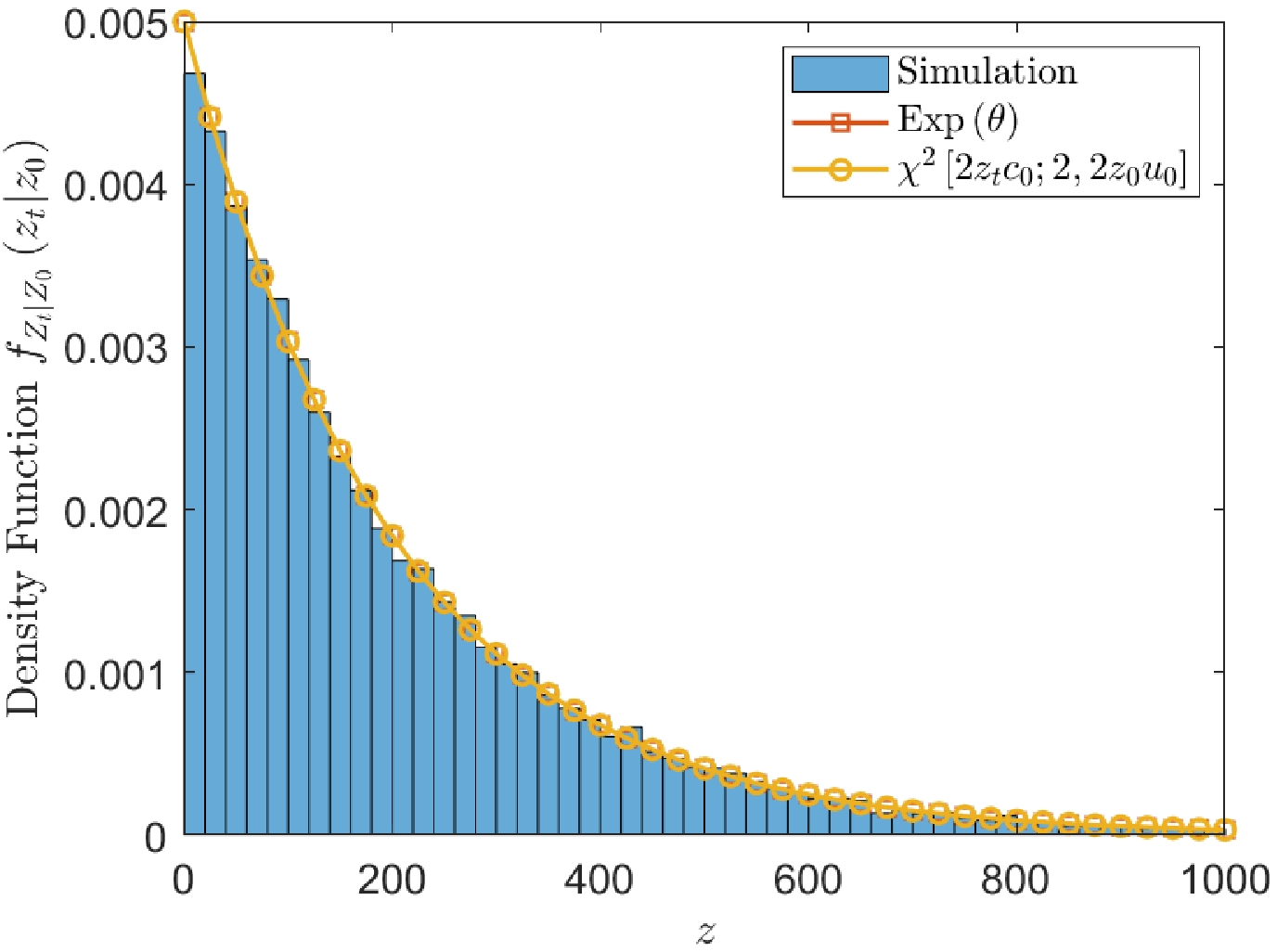}
\caption{The conditional probability density of the squared distance $f_{Z_t|Z_0}\left(z_t|z_0\right)$ at time $t=10$s; mobility parameters $\tau=1$s and $D=100$m$^2/$s.}
\label{fig:cond-pdf-dist-10}
\end{figure}

Starting from the mobility model described in section~\ref{subsec:mobility}, we derive the SDE of $Z_t$ stated in the following proposition. 
\begin{prop}
The stochastic differential equation of the squared distance $Z_t$ between two nodes $1$ and $2$ moving randomly in a two-dimensional plane according to independent identical OU processes with mobility parameters $\tau$ and $D$ is 
\begin{equation}\label{eq:sde_z}
\mathrm{d}Z_t=k\left(\theta-Z_t\right)\mathrm{d}t+\sigma\sqrt{Z_t}\mathrm{d}W_t, \quad Z_{0}=z_0,
\end{equation}
where $k=2/\tau$, $\theta = 2D\tau$, $\sigma=2\sqrt{2D}$, $W_t$ is a standard Brownian motion, and $z_0\geq 0$ is the starting point of the process.
\end{prop}
\begin{IEEEproof}
See Appendix~\ref{Appen-Dist-SDE}.
\end{IEEEproof}
%%%%%%%%%%%%%%%%%%%%%%%%%%%%%%%%%%%%%%%%
%%%%%%%%%%%%%%%%%%%%%%%%%%%%%%%%%%%%%%%%
\begin{figure}[t]
\centering
\includegraphics[width = \columnwidth]{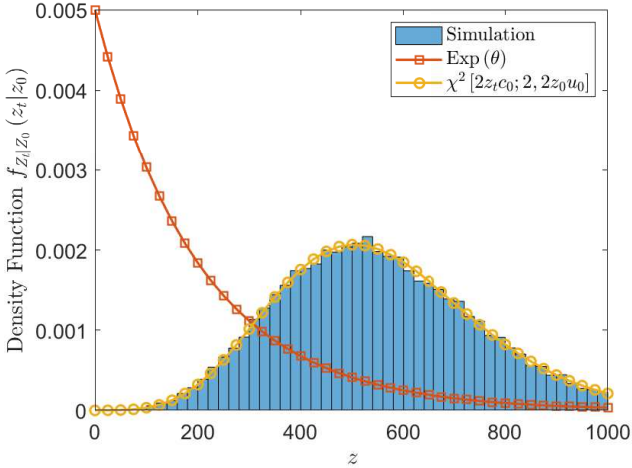}
\caption{The conditional probability density of the squared distance $f_{Z_t|Z_0}\left(z_t|z_0\right)$ at time $t=0.2$s; mobility parameters $\tau=1$s and $D=100$m$^2/$s.}
\label{fig:cond-pdf-dist-2}
\end{figure}
%%%%%%%%%%%%%%%%%%%%%%%%%%%%%%%%%%%%%%%%
%%%%%%%%%%%%%%%%%%%%%%%%%%%%%%%%%%%%%%%%
The stochastic process $Z_t$ solving~\eqref{eq:sde_z} belongs to the Cox-Ingersoll-Ross (CIR) family of diffusions. They were first introduced in finance to model short-term interest rates~\cite{cox1985theory}. By performing particular space-time changes, $Z_t$ can be represented as a Bessel-squared process, $BESQ_{z_0}^{\delta}$, with dimension $\delta=2$. The transition densities of Bessel-squared processes are known explicitly~\cite{jeanblanc2009mathematical}; therefore, the transition density of $Z_t$ can be determined precisely as follows.
\begin{cor}
The squared distance at time $s$, conditioned on its value at the current time $t$, follows a non-central chi-square distribution, $\chi^2\left[2z_sc;2,2z_tu\right]$, with $2$ degrees of freedom and parameter of non-centrality $2z_tu$. Its transition pdf can be expressed as 
\begin{equation}\label{eq:nonchi}
f_{Z_s|Z_t}\left(z_s|z_t\right)=ce^{-(z_sc+z_tu)}\mathrm{I}_0\left(2\sqrt{z_tuz_sc}\right),
\end{equation}
where $c=\frac{1}{\theta\left(1-e^{-k(s-t)}\right)}$, $u=ce^{-k(s-t)}$, and $\mathrm{I}_0$ is the modified Bessel function of the first kind with order zero.
\end{cor}
\begin{IEEEproof}
The result follows from a time-space transformation of the transition density of the corresponding Bessel-squared process~\cite{jeanblanc2009mathematical}.
\end{IEEEproof} 
The parameters $k$ and $\theta$ influence the behavior of the process $Z_t$ in several ways~\cite{going2003survey}. First, if $k,\theta>0$, then~\eqref{eq:sde_z} admits a unique solution. Second, if $2k\theta\geq\sigma^2$, the stochastic process $Z_t$ is strictly positive for $t>0$, and never hits zero~\cite{cox1985theory, feller1951two}. These conditions are always verified in our model for every value of $\tau$ and $D$. Therefore, $Z_t$ is a mean-reverting diffusion process with speed of adjustment $k$ and long-term average $\theta$. Next, we make two remarks. 
\begin{rem}\label{rem:limiting-distrib}
If $k,\theta>0$, then as $s\rightarrow\infty$ the conditional density will approach an exponential distribution, $f_{Z_{\infty}|Z_t}\left(z\right)\sim\mathsf{Exp}\left(\theta\right)$, with mean and variance equal to $\theta$ and $\theta^2$, respectively. The result follows by taking the limit for $s\rightarrow\infty$ of~\eqref{eq:nonchi}.
\end{rem}
\begin{rem}\label{rem:limiting-distrib}
The stochastic process $Z_t$ possesses the Markov property. The result follows from the fact that $Z_t$ is a Bessel-squared process, and Bessel (squared) processes are Markov processes~\cite{jeanblanc2009mathematical}.
\end{rem}
We verify the validity of our theoretical analysis by running different simulations. In Fig.~\ref{fig:cond-pdf-dist-2} we plot the conditional probability density of the squared distance $f_{Z_t|Z_0}\left(z_t|z_0\right)$ at time $t=0.2$s, whereas, in Fig.~\ref{fig:cond-pdf-dist-10} we plot the same function but at time $t=10$s. To generate these plots, we fix the starting point of the squared distance process, $Z_0$, and then perform Monte Carlo simulations to estimate the conditional distribution function empirically. At the beginning of the simulations, $t=0.2$s, we observe that the conditional distribution follows a non-central chi-square distribution. Eventually, as time goes by, it approaches the limiting distribution.
%%%%%%%%%%%%%%%%%%%%%%%%%%%%%%%%%%%%%%%%
%%%%%%%%%%%%%%%%%%%%%%%%%%%%%%%%%%%%%%%%
%%%%%%%%%%%%%%%%%%%%%%%%%%%%%%%%%%%%%%%%
%%%%%%%%%%%%%%%%%%%%%%%%%%%%%%%%%%%%%%%% 

\begin{figure}[t]
\centering
\includegraphics[width = \columnwidth]{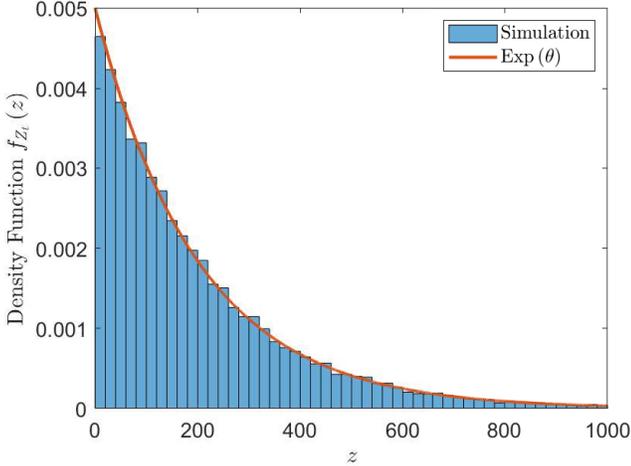}
\caption{The probability density, $f_{Z_t}\left(z\right)$, of the squared distance at any time $t$; mobility parameters $\tau=1$s and $D=100$m$^2/$s.}
\label{fig:pdf-dist-stat}
\end{figure}

Now, if the initial condition of the process, $Z_0$, is drawn according to the  limiting distribution, then the process $\{Z_t, t\geq0\}$ is stationary. Beginning with~\eqref{eq:nonchi} and averaging over $Z_0$, the density function of the squared distance can be evaluated to yield the result stated in the following proposition.
\begin{prop}\label{prop:1}
At every time $t$, when the initial condition of the process $Z_0$ is drawn according to the limiting distribution, the squared distanced between the nodes $1$ and $2$, $Z_t$, is exponentially distributed. Its pdf can be written as
\begin{equation}\label{eq:sta_distrib_dist}
f_{Z_t}\left(z\right)=\frac{1}{\theta}e^{-z/\theta},
\end{equation}
where $\theta = 2D\tau$. 
\end{prop}
\begin{IEEEproof}
See Appendix~\ref{Appen-Dist-pdf}.
\end{IEEEproof}
Equation~\eqref{eq:sta_distrib_dist} shows that $\{Z_t, t\geq0\}$ is a stationary stochastic process, since its distribution does not depend on time $t$. Therefore, we can conclude that the squared distance between two nodes $1$ and $2$ forms a stationary Markov process with transition probability density given by~\eqref{eq:nonchi}, and steady-state pdf given by~\eqref{eq:sta_distrib_dist}. To the best of our knowledge, this is the most complete mathematical analysis reported in the literature proving that the distance process inherits the stationary and Markov properties of the OU process modeling node displacements along the $x$ and $y$ coordinates. As an example, in Fig.~\ref{fig:pdf-dist-stat}, we plot the pdf~\eqref{eq:sta_distrib_dist}. It should be noted that in this figure the starting point of the process, $Z_0$, is drawn independently from the limiting distribution in each simulation trial. We observe an excellent agreement between the theoretical result and the Monte Carlo estimate.
%%%%%%%%%%%%%%%%%%%%%%%%%%%%%%%%%%%%%%%%
%%%%%%%%%%%%%%%%%%%%%%%%%%%%%%%%%%%%%%%%
%%%%%%%%%%%%%%%%%%%%%%%%%%%%%%%%%%%%%%%%
%%%%%%%%%%%%%%%%%%%%%%%%%%%%%%%%%%%%%%%%
\subsection{Autocovariance Function}
\begin{figure}[t]
\centering
\includegraphics[width = \columnwidth]{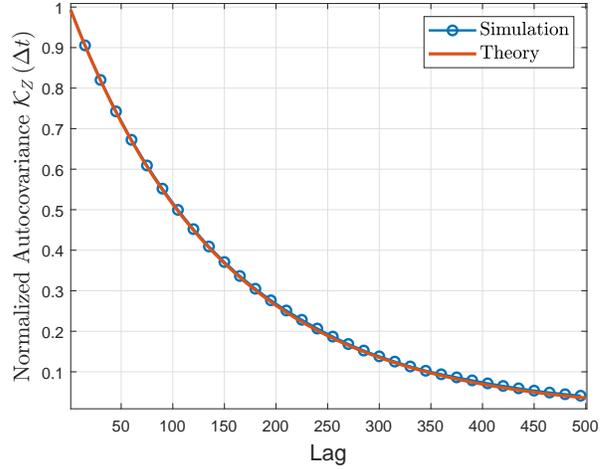}
\caption{The normalized autocovariance function $\mathcal{K}_{Z}\left(\Delta t\right)$; mobility parameters $\tau=0.1$s and $D=100$m$^2/$s.}
\label{fig:cov-dist}
\end{figure}

Another statistical property of interest is the stationary autocovariance function of the process $\{Z_t, t\geq0\}$. It gives a measure of dependency of the random process to its delayed version, as a function of the time-lag, in the steady-state. Given that $Z_t$ is stationary, this function depends only on the time shift $s-t=\Delta t$.
\begin{prop}\label{prop:autocov}
The stationary autocovariance function of the process $\{Z_t, t\geq0\}$ is given by
\begin{align}
K_{Z}\left(\Delta t\right)=\theta^2e^{-k\Delta t},
\end{align}
where $k=2/\tau$, and $\theta = 2D\tau$.
\end{prop}
\begin{IEEEproof}
See Appendix~\ref{Appen-Dist-acf}.
 \end{IEEEproof}
This result implies that the autocovariance function tends toward zero as the time-lag increases. This is perfectly intuitive since an increase in time-lag would yield a decrease in the correlation of the random variables $Z_s$ and $Z_t$. In Fig.~\ref{fig:cov-dist} we plot the normalized autocovariance function $K_{Z}\left(\Delta t\right)$ versus the time-lag, and we observe how the autocovariance function decays exponentially with the time-lag.
\begin{rem}
The stationary autocovariance function $K_{Z}\left(\Delta t\right)$ is proportional to the autocovariance function of the OU process $K_{S}\left(\Delta t\right)$ in~\eqref{eq:autocorr}. 
\end{rem}
%%%%%%%%%%%%%%%%%%%%%%%%%%%%%%%%%%%%%%%%
%%%%%%%%%%%%%%%%%%%%%%%%%%%%%%%%%%%%%%%%
%%%%%%%%%%%%%%%%%%%%%%%%%%%%%%%%%%%%%%%%
%%%%%%%%%%%%%%%%%%%%%%%%%%%%%%%%%%%%%%%%
%%%%%%%%%%%%%%%%%%%%%%%%%%%%%%%%%%%%%%%%
%%%%%%%%%%%%%%%%%%%%%%%%%%%%%%%%%%%%%%%%
\subsection{Distance Process}
\begin{figure}[t]
\centering
\includegraphics[width = \columnwidth]{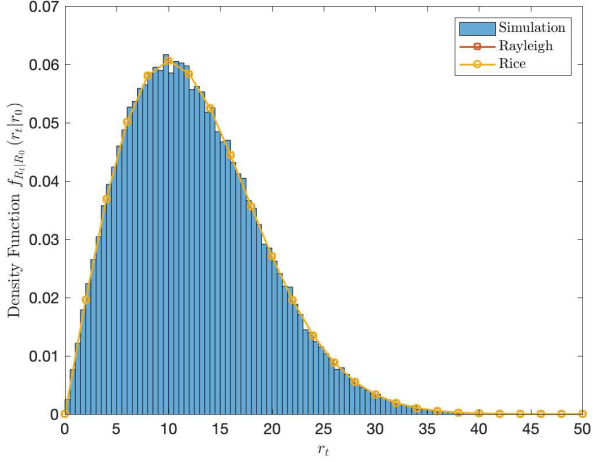}
\caption{The conditional probability density of the distance $f_{R_t|R_0}\left(r_t|r_0\right)$ at time $t=10$s; mobility parameters $\tau=1$s and $D=100$m$^2/$s.}
\label{fig:cond-pdf-r-dist-10}
\end{figure}

Here, we analyze the distance between nodes, $R_t=\sqrt{Z_t}$. From the theory of Bessel processes~\cite{jeanblanc2009mathematical}, the square root of $BESQ_{z_0}^{\delta}$ is also a Bessel process of dimension $\delta$ and starting point $r_0=\sqrt{z_0}$.
\begin{cor}
The transition probability density of the distance process $R_t$ follows a Rice distribution, while its steady-state is Rayleigh distributed, i.e.,
 \begin{equation}
f_{R_s|R_t}\left(r_s|r_t\right)=\frac{r_s}{b^2}\exp\left[\frac{-\left(r_s^2+r_t^2a^2\right)}{2b^2}\right]\mathrm{I}_0\left(\frac{r_sr_ta}{b^2}\right),
\end{equation}
and
\begin{equation}
f_{R_t}\left(r\right)=\frac{2r}{\theta}e^{-r^2/\theta}
\end{equation}
respectively, with $\theta = 2D\tau$, $a^2=u/c$, and $b^2=1/2c$.
\end{cor}
\begin{IEEEproof}
Apply the transformation of random variables $R_t=\sqrt{Z_t}$ and the result follows.
\end{IEEEproof} 
In Fig.~\ref{fig:cond-pdf-r-dist-10} we plot the conditional probability density of the distance $f_{R_t|R_0}\left(r_t|r_0\right)$ for the particular time instance $t=10$s. In a similar fashion as for the squared distance, to generate this plot we fix the starting point $R_0$ and then perform Monte Carlo simulations to estimate the conditional distribution function empirically. As expected, we observe the same behavior, i.e., at the beginning of the simulations, the conditional distribution follows a Rice distribution. Eventually, as time goes by, it approaches the steady-state which is Rayleigh distributed.
To conclude, in a system where nodes move randomly in a two-dimensional plane according to an OU process both the distance $R_t$ and the squared distance $Z_t$ form stationary Markov processes.
 %%%%%%%%%%%%%%%%%%%%%%%%%%%%%%%%%%%%%%%%
%%%%%%%%%%%%%%%%%%%%%%%%%%%%%%%%%%%%%%%%
%%%%%%%%%%%%%%%%%%%%%%%%%%%%%%%%%%%%%%%%
%%%%%%%%%%%%%%%%%%%%%%%%%%%%%%%%%%%%%%%%
%%%%%%%%%%%%%%%%%%%%%%%%%%%%%%%%%%%%%%%%
%%%%%%%%%%%%%%%%%%%%%%%%%%%%%%%%%%%%%%%%
\subsection{SNR process}
Now, starting from the SDE of $Z_t$ in~\eqref{eq:sde_z} we will derive the SDE of $N_t$ by applying It\^{o}'s formula, which enables us to find the differential of a scalar function of a stochastic process~\cite{oksendal2003stochastic}. This formula is valid when the function is $\mathcal{C}^2$, i.e., it is twice differentiable. The SNR $N_t$ is not $\mathcal{C}^2$ at the origin; but, in our system, $Z_t$ never reaches zero for $t>0$, and $Z_0$ is drawn according to the limiting distribution, i.e., $Z_0\sim\mathsf{Exp}\left(\theta\right)$. Therefore, since $Z_t$ never hits the origin for $t\geq0$, we can still apply It\^{o}'s formula and obtain the following result.
\begin{prop}
The stochastic differential equation of the link SNR $N_t$ in mobile networks where nodes experience an OU mobility model is 
\begin{align}\label{eq:sdeSNR}\nonumber
\mathrm{d}N_t=&\left[\frac{k\eta}{2}N_t-\frac{k\theta\eta}{2}N_t^{1+2/\eta}+\frac{\sigma^2}{2}\frac{\eta}{2}\left(1+\frac{\eta}{2}\right)N_t^{1+2/\eta}\right]\mathrm{d}t\\
&- \frac{\sigma\eta}{2}N_t^{1+1/\eta}\mathrm{d}W_t, \quad N_{0}=Z_0^{-\eta/2},
\end{align}
where $k=2/\tau$, $\theta = 2D\tau$, $\sigma=2\sqrt{2D}$, and $W_t$ is a standard Brownian motion.
\end{prop}
\begin{IEEEproof}
See Appendix~\ref{Appen-SNR-SDE}.
\end{IEEEproof} 
Furthermore, the pdf and the cdf of the instantaneous SNR $N_t$ can be determined precisely as follows. 
\begin{cor}
The density function of the link SNR is
\begin{equation}\label{eq:sta_distrib_snr}
f_{N_t}\left(\rho\right)=\frac{2}{\eta\theta}\rho^{-2/\eta - 1}e^{-\frac{\rho^{-2/\eta}}{\theta}},\quad \rho>0,
\end{equation}
and its cumulative distribution is
\begin{equation}\label{eq:cdf_snr}
F_{N_t}\left(\rho\right)=e^{-\frac{\rho^{-2/\eta}}{\theta}},
\end{equation}
where $\theta = 2D\tau$. 
\end{cor}
\begin{IEEEproof}
To obtain the pdf, apply the principle of random variable transformation to~\eqref{eq:snr_dist} and the result follows. The cdf, instead, is obtained by integrating the pdf in~\eqref{eq:sta_distrib_snr}. 
\end{IEEEproof} 
%%%%%%%%%%%%%%%%%%%%%%%%%%%%%%%%%%%%%%%%
%%%%%%%%%%%%%%%%%%%%%%%%%%%%%%%%%%%%%%%%
\begin{figure}[t]
\centering
\includegraphics[width = \columnwidth]{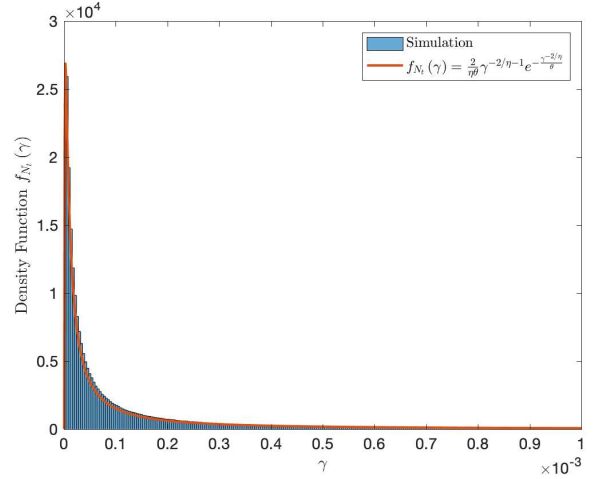}
\caption{The probability density of the link SNR $N_t$ for $\eta=4$; mobility parameters $\tau=1$s and $D=100$m$^2/$s.}
\label{fig:snr}
\end{figure}
%%%%%%%%%%%%%%%%%%%%%%%%%%%%%%%%%%%%%%%%
%%%%%%%%%%%%%%%%%%%%%%%%%%%%%%%%%%%%%%%%
To check the accuracy of the derived pdf in~\eqref{eq:sta_distrib_snr}, in Fig.~\ref{fig:snr} we plot the pdf of the link SNR for $\eta=4$. We observe an excellent match between the theoretical result and the Monte Carlo simulate. It is important to note that the mean and the variance of the link SNR are undefined because the integral does not converge. 

\begin{figure}[t]
\centering
\includegraphics[width = \columnwidth]{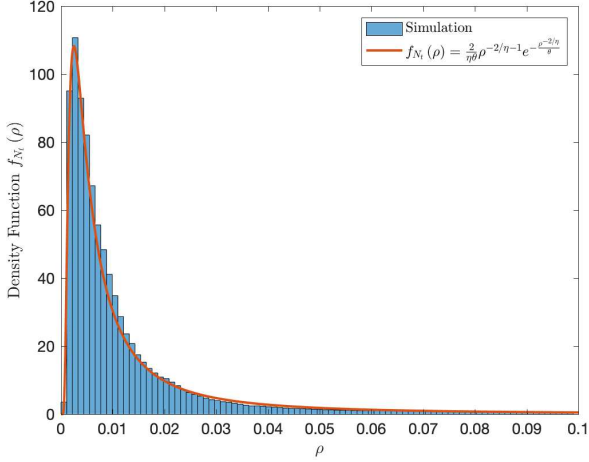}
\caption{The probability density of the link SNR $N_t$ for $\eta=2$; mobility parameters $\tau=1$s and $D=100$m$^2/$s.}
\label{fig:sde-snr}
\end{figure}

Eq.~\eqref{eq:sdeSNR} provide a useful way to simulate values of the SNR $N_t$. In particular, we can simulate a discretized process $\{N_{k}, k \in \mathbb{N}\}$ for any time step $k=t_0 + k\Delta t$, $k\in \mathbb{N}$ and $\Delta t > 0$. The smaller the value of $\Delta t$, the closer our discretized path will be to the continuous-time path of~\eqref{eq:sdeSNR}. We perform numerical simulations to check the accuracy of the derived SDE. While there are a number of discretization schemes available to simulate the SDE in~\eqref{eq:sdeSNR}, we used the simplest and most common scheme, the Euler scheme. In Fig.~\ref{fig:sde-snr}, we plot the density of the discretized path of the link SNR and the pdf in~\eqref{eq:sta_distrib_snr} when $\eta=2$. 

Understanding whenever the link SNR $N_t$ is greater than some certain threshold is very important for connectivity issues. It determines if a successful communication can be established between the two nodes. From~\eqref{eq:cdf_snr}, the probability that two nodes $1$ and $2$ are connected is given by
\begin{equation}
\mathbb{P}\left(N_{t}\geq \rho_{th}\right)=1-e^{-\frac{\rho_{th}^{-2/\eta}}{\theta}}. 
\end{equation}
From a system design perspective, this is a very important result, because it provides an explicit relation between the connectivity probability and various system parameters, such as the transmit power, wavelength, bandwidth, and the noise spectral density. 
%%%%%%%%%%%%%%%%%%%%%%%%%%%%%%%%%%%%%%%%
%%%%%%%%%%%%%%%%%%%%%%%%%%%%%%%%%%%%%%%%
\subsection{Bivariate Distribution of the Link SNR}
In many situations, we may be interested in the future states of the link SNR. For instance, we may ask what is the probability that the link SNR is lower than some certain threshold for two different time instances. To that end, we derive the joint cdf of the link SNR random variables $N_s$ and $N_t$, with $s>t$, as follows.
\begin{prop}
The bivariate cumulative distribution of the link SNR random variables $N_s$ and $N_t$, $s>t$, is
\begin{align}\label{eq:joint_cdf_snr}
F_{N_s,N_t}\left(\rho_s,\rho_t\right)=&\frac{1}{\theta}\sum_{j=0}^{\infty}\frac{u^j}{c^{j+1}j!\Gamma\left(j+1\right)}\\\nonumber
&\times\gamma\left(j+1,c\rho_s^{-2/\eta}\right)\gamma\left(j+1,c\rho_t^{-2/\eta}\right),
\end{align}
where $c=\frac{1}{\theta\left(1-e^{-k(s-t)}\right)}$, $u=ce^{-k(s-t)}$, $\Gamma\left(\cdot\right)$ is the gamma function, and $\gamma\left(a,x\right)=\int_0^x v^{a-1}e^{-v}\mathrm{d}v$ is the lower incomplete gamma function.
\end{prop}
\begin{IEEEproof}
See Appendix~\ref{Appen-joint-cdf-snr}.
\end{IEEEproof} 
\begin{cor}
The bivariate density function of the link SNR random variables $N_s$ and $N_t$, $s>t$, is
\begin{align}\label{eq:joint_pdf_snr}
f_{N_s,N_t}\left(\rho_s,\rho_t\right)=&\frac{4}{\eta^2\theta}\sum_{j=0}^{\infty}\frac{u^jc^{j+1}}{j!\Gamma\left(j+1\right)}\\\nonumber
&\times\left[\left(\rho_s\rho_t\right)^{-2/\eta}\right]^je^{-c\left(\rho_s^{-2/\eta}+\rho_t^{-2/\eta}\right)},
\end{align}
where $c=\frac{1}{\theta\left(1-e^{-k(s-t)}\right)}$, $u=ce^{-k(s-t)}$, and $\Gamma\left(\cdot\right)$ is the gamma function.
\end{cor}
\begin{IEEEproof}
The result follows from differentiating the joint cumulative distribution function~\eqref{eq:joint_cdf_snr} with respect to $\rho_s$ and $\rho_t$, i.e., $f_{N_s,N_t}\left(\rho_s,\rho_t\right)=\frac{\delta^2F_{N_s,N_t}\left(\rho_s,\rho_t\right)}{\delta\rho_s\delta\rho_t}$.
\end{IEEEproof} 
%%%%%%%%%%%%%%%%%%%%%%%%%%%%%%%%%%%%%%%%
%%%%%%%%%%%%%%%%%%%%%%%%%%%%%%%%%%%%%%%%
%%%%%%%%%%%%%%%%%%%%%%%%%%%%%%%%%%%%%%%%
%%%%%%%%%%%%%%%%%%%%%%%%%%%%%%%%%%%%%%%%
%%%%%%%%%%%%%%%%%%%%%%%%%%%%%%%%%%%%%%%%
%%%%%%%%%%%%%%%%%%%%%%%%%%%%%%%%%%%%%%%%
\section{Statistical Properties of Mobile Networks subject to Rayleigh Fading}\label{sec:with-fading}
In this section, we generalize our analysis to take into account small-scale fading modeled by a Rayleigh random variable.  Through this extension, the distribution of the link SNR depends both on the squared distance $Z_t$, which is governed by the mobility model, and the environmental factors controlling the channel between devices captured by $G_t$. It is clear from~\eqref{eq:inst_snr} that the link SNR has a compound  probability distribution, i.e., $ N_{t}\sim\mathsf{Exp}\left(\Upsilon_t\right)$, where $\Upsilon_t=Z_{t}^{\eta/2}$. On that account, the unconditional density function of the link SNR can be evaluated to yield the result stated in the following proposition.
\begin{prop}\label{prop:dist_snr}
In a system where nodes move randomly according to an OU process and Rayleigh fading affects their connections, the pdf of the link SNR at any time $t$ for rational path loss exponents $\eta$ is given by
\begin{multline}\label{eq:pdf-snr-eta}
f_{N_t}\left(\rho\right)=\frac{2q}{p\theta}\frac{p^{\frac{3}{2}+\frac{p}{2q}}}{\sqrt{2q}}\left(2\pi\right)^{1-q-\frac{p}{2}}\theta^{1+\frac{p}{2q}}\\
\times\MeijerG*{2q}{p}{p}{2q}{\frac{2pq-2q-p}{2pq}, \dots, \frac{-p}{2pq}}{0, \dots, \frac{2q-1}{2q}}{\left(\frac{\rho}{2q}\right)^{2q}\left(p\theta\right)^{p}},\quad  \rho\geq 0,
\end{multline}
where $\MeijerG*{m}{n}{s}{t}{u_1, \dots, u_s}{v_1, \dots, v_t}{z}$ is the Meijer $G$ furcation, $\theta=2D\tau$, $p,q\in \mathbb{Z}^{+}$ so that $\eta=p/q$ is a positive rational number. 
\end{prop}
\begin{IEEEproof}
See Appendix~\ref{Appen-MT}.
\end{IEEEproof}
%%%%%%%%%%%%%%%%%%%%%%%%%%%%%%%%%%%%%%%%
%%%%%%%%%%%%%%%%%%%%%%%%%%%%%%%%%%%%%%%%
\begin{figure}[t]
\centering
\includegraphics[width = \columnwidth]{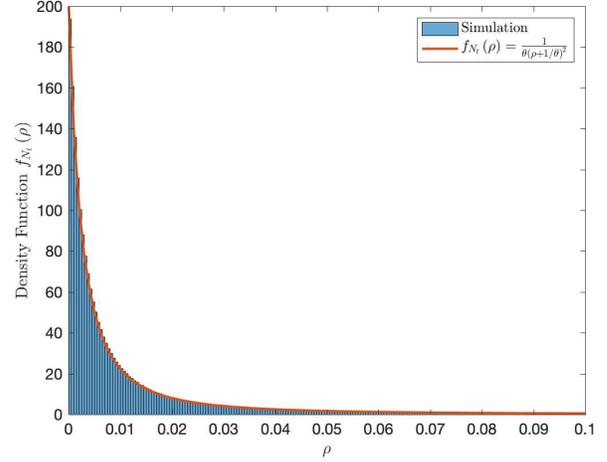}
\caption{The probability density of the link SNR $N_t$ for $\eta=2$; mobility parameters $\tau=1$s and $D=100$m$^2/$s; channel parameters $\nu_{max}=100$Hz, sampling rate=$0.0003$s.}
\label{fig:snr_fading}
\end{figure}
%%%%%%%%%%%%%%%%%%%%%%%%%%%%%%%%%%%%%%%%
%%%%%%%%%%%%%%%%%%%%%%%%%%%%%%%%%%%%%%%%

Eq.~\eqref{eq:pdf-snr-eta} provides an explicit relation between the pdf of the link SNR, mobility parameters $\tau$, $\theta$, and path loss exponent $\eta$. A number of interesting points can be noted from this expression. First, it indicates that the link SNR $\{N_t, t\geq0\}$ is a first-order stationary stochastic process, since its distribution does not depend on time $t$. Second, this is the most complete mathematical analysis of the distribution of the link SNR reported in the literature for a system subject to Rayleigh fading and OU Mobility accounting for rational path loss exponent. Indeed, since the path loss exponent is an experimentally estimated parameter, it is, by definition, rational in practice due to finite precision measurement equipment. Although the link SNR distribution is given in terms of the Meijer $G$ function, it can be easily evaluated using numerical software such as Mathematica for any given inputs. It should be noted that~\eqref{eq:pdf-snr-eta} reduces to the following expressions for the special cases $\eta=2$
\begin{equation}\label{eq:pdf_snr_eta2}
f_{N_t}\left(\rho\right)=\frac{1}{\theta\left( \rho+1/\theta\right)^{2}},\quad  \rho\geq 0,
\end{equation}
and $\eta=4$
\begin{equation}\label{eq:pdf_snr_eta4}
f_{N_t}\left(\rho\right)=\frac{\sqrt{\pi}e^{1/4\rho\theta^2}\left(1+2\theta^2\rho\right)\mathrm{Erfc}\left(\frac{1}{2\theta\sqrt{\rho}}\right)-2\theta\sqrt{\rho}}{4\theta^2\rho^{5/2}}
\end{equation}
where $\rho\geq0$, and $\mathrm{Erfc}\left(x\right)=\frac{2}{\sqrt{\pi}}\int_x^{\infty}e^{-t^2}\mathrm{d}t$ is the complementary error function. However, for any other value of $\eta$, the expression given in proposition~\ref{prop:dist_snr} is the most compact, accessible form. Eq.~\eqref{eq:pdf_snr_eta2} corresponds to a shifted Pareto distribution with shape parameter $1$ and scale parameter $1/\theta$. This is a heavy-tail distribution, with undefined mean and variance given that the shape parameter is equal to one. 

\begin{figure}[t]
\centering
\includegraphics[width = \columnwidth]{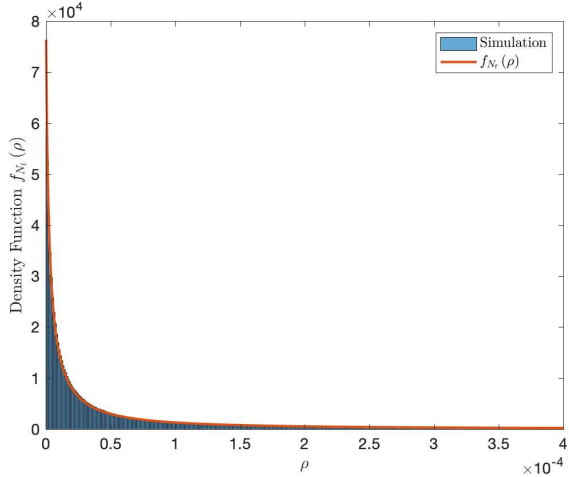}
\caption{The probability density of the link SNR $N_t$ for $\eta=4$; mobility parameters $\tau=1$s and $D=100$m$^2/$s; channel parameters $\nu_{max}=100$Hz, sampling rate=$0.0003$s.}
\label{fig:snr_fading_4}
\end{figure}

We perform numerical simulations to check the accuracy of the derived distribution formulae in~\eqref{eq:pdf_snr_eta2} and in~\eqref{eq:pdf_snr_eta4}. To simulate the Rayleigh fading channel, we refer to the autoregressive stochastic model presented in~\cite{baddour2005autoregressive}. In Fig.~\ref{fig:snr_fading}  and Fig.~\ref{fig:snr_fading_4} we plot the pdf of the link SNR when $\eta=2$ and $\eta=4$, respectiovely. We observe an excellent agreement between the Monte Carlo result and the corresponding one from the mathematical analysis.

For the particular case $\eta=2$, the cumulative distribution of the link SNR can be evaluated in closed-form by integrating~\eqref{eq:pdf_snr_eta2}
\begin{equation}
F_{N_t}\left(\rho\right)=\frac{\rho\theta}{1+\rho\theta},
\end{equation}
Then, the probability that two nodes $1$ and $2$ are connected is given by
\begin{equation}
\mathbb{P}\left(N_t\geq \rho_{th}\right)=\frac{1}{1+\rho_{th}\theta},
\end{equation}
with $\theta=2D\tau$. This is a significant result because it provides some insight to the relationship between the link connectivity probability, mobility parameters $\tau$ and $\theta$, and the SNR threshold $\rho_{th}$, which depends on the communication hardware, as well as the modulation and coding scheme of the wireless system. 
%%%%%%%%%%%%%%%%%%%%%%%%%%%%%%%%%%%%%%%%
%%%%%%%%%%%%%%%%%%%%%%%%%%%%%%%%%%%%%%%%
%%%%%%%%%%%%%%%%%%%%%%%%%%%%%%%%%%%%%%%%
%%%%%%%%%%%%%%%%%%%%%%%%%%%%%%%%%%%%%%%%
\section{Discussion and Conclusions}\label{sec:conclusions}
In this work, we derived closed-form expressions for the statical properties of the link SNR and the separation distance in systems subject to Rayleigh Fading and Ornstein-Uhlenbeck Mobility. We started our analysis by first considering the case when there is no signal fading affecting the link between nodes. In this scenario, the statistical properties of the link SNR are entirely determined by the squared distance between the two nodes. As main contributions, we provided a full statistical description of the squared distance process, including its distribution and autocorrelation function, and showed that it forms a stationary Markov process. Then, we derived closed-form expressions for the pdf, the cdf, the bivariate pdf, and the bivariate cdf of the link SNR. Next, we extended our analysis to take into account variations in the propagation channel (e.g., fading), and calculated the pdf of the link SNR for rational path loss exponents $\eta$. To the best of our knowledge, this is the most complete mathematical analysis reported in the literature for the distribution of the link SNR in mobile wireless systems. We also computed expressions for the connectivity probability in closed form for both fading and non-fading scenarios. Finally, we performed extensive simulations to check the accuracy of the proposed mathematical analysis.

\begin{figure}[t]
\centering
\includegraphics[width = \columnwidth]{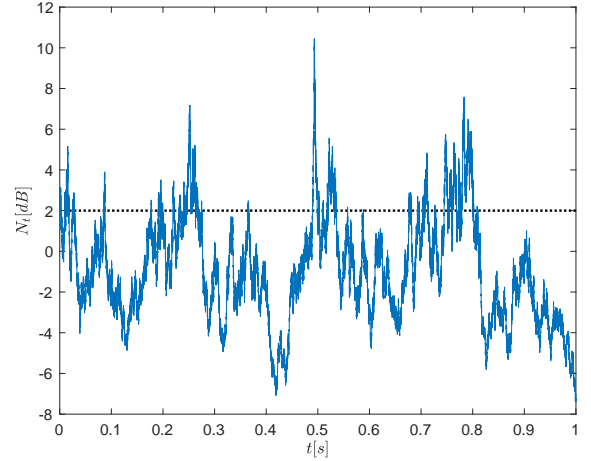}
\caption{Simulation of a SNR path for $\eta=2$ in mobile networks without fading; $\rho_{th} = 2$dB; mobility parameters $\tau=0.6$s and $D=4$m$^2/$s.}
\label{fig:snr_path}
\end{figure}

Characterizing and managing the SNR variations users would see across mobile ad hoc networks is a challenging but important problem towards understanding network stability and connectivity. Therefore, the results derived in this paper could be helpful to quantify the coverage and outage durations that each user will experience in the network. For instance, given an SNR threshold, as shown in Fig.~\ref{fig:snr_path}, one can characterize the temporal characteristics of the on/off level crossing process associated with the SNR being above and below the threshold. 

%%%%%%%%%%%%%%%%%%%%%%%%%%%%%%%%%%%%%%%%
%%%%%%%%%%%%%%%%%%%%%%%%%%%%%%%%%%%%%%%%
%%%%%%%%%%%%%%%%%%%%%%%%%%%%%%%%%%%%%%%%
%%%%%%%%%%%%%%%%%%%%%%%%%%%%%%%%%%%%%%%%

Finally, our results provide useful insight and analytical tools that can be used to develop a framework to evaluate link stability in systems subject to Rayleigh Fading and Ornstein-Uhlenbeck Mobility. These systems pose several design issues due to the dynamic characteristics of their underlying topology. Particularly, the routing protocol faces strong challenges as connections between nodes are established and broken intermittently. It is, therefore, imperative to quantify topological uncertainty when designing and implementing these systems in any real-world application. In our previous work~\cite{cika2019quantifying}, we presented a mobility metric for the evaluation of the link stability in mobile ad hoc networks in the absence of signal fading. Our next goal is to generalize this stability metric to take into account the effects of the physical characteristics of the underlying wireless channel. 
%%%%%%%%%%%%%%%%%%%%%%%%%%%%%%%%%%%%%%%%
%%%%%%%%%%%%%%%%%%%%%%%%%%%%%%%%%%%%%%%%
\appendices
 \section{Proof of the stochastic differential equation of the squared distance process}
\label{Appen-Dist-SDE}
From the mobility model described in section~\ref{subsec:mobility}, $\{X_t, t\geq0\}$ and $\{Y_t, t\geq0\}$ are two identical independent OU processes. Using It\^{o}'s formula, $\left(\right.$Theorem $4.1.2$ in~\cite{oksendal2003stochastic}$\left.\right)$, we can compute the SDEs of the squared stochastic processes as follows
\begin{align}\label{eq:sde_x^2}
\mathrm{d}X_t^2&=\left(2D-\frac{2}{\tau}X_t^2\right)\mathrm{d}t+2\sqrt{2D}X_t\mathrm{d}W^x_t,\\\label{eq:sde_y^2}
\mathrm{d}Y_t^2&=\left(2D-\frac{2}{\tau}Y_t^2\right)\mathrm{d}t+2\sqrt{2D}Y_t\mathrm{d}W^y_t,
\end{align}
where $\{W^x_t, t\geq 0\}$ and $\{W^y_t, t\geq 0\}$ are two independent standard Brownian motions.

Now, let
\begin{equation}
B_t=\int_0^t X_u\mathrm{d}W^x_u+Y_u\mathrm{d}W^y_u.
\end{equation}
$B_t$ is a stochastic integral with respect to a Brownian motion with quadratic variation given by
\begin{equation}
\langle B \rangle_t=\int_0^t\left[ \left(X_u\right)^2+\left(Y_u\right)^2\right]\mathrm{d}u=\int_0^t Z_u \mathrm{d}u.
\end{equation}
Consequently, by Levy's characterization theorem $\left(\right.$Theorem $8.6.1$ in~\cite{oksendal2003stochastic}$\left.\right)$, the stochastic process
\begin{equation}
W_t=\int_0^t\frac{1}{\sqrt{Z_u}} \left(X_u\mathrm{d}W^x_u+Y_u\mathrm{d}W^y_u\right)
\end{equation}
is a standard Brownian motion. This enables to write the SDE of $Z_t$ as
\begin{equation}\label{ap:sdez}
\mathrm{d}Z_t=k\left(\theta-Z_t\right)\mathrm{d}t+\sigma\sqrt{Z_t}\mathrm{d}W_t, \quad Z_{0}=z_0,
\end{equation}
where $k=2/\tau$, $\theta = 2D\tau$, $\sigma=2\sqrt{2D}$, and $Z_0$ is the starting point of the process. The solution to this SDE exists, and is unique if $k,\theta > 0$. It is given by
\begin{multline}\label{eq:SDEsol}
Z_t=z_0e^{-k t}+\theta\left(1-e^{-k t}\right)\\
+\sigma\int_0^t\sqrt{Z_u}e^{k\left(u-t\right)}\mathrm{d}W_u.
\end{multline}
This expression is obtained by multiplying both sides of~\eqref{ap:sdez} by $e^{k t}$ and then integrating over the time interval $[0,t]$.
%%%%%%%%%%%%%%%%%%%%%%%%%%%%%%%%%%%%%%%%
%%%%%%%%%%%%%%%%%%%%%%%%%%%%%%%%%%%%%%%%
\section{Proof of the density function of the squared distance process}
\label{Appen-Dist-pdf}
Here, we provide a brief outline for the derivation of the pdf expression in~\eqref{eq:sta_distrib_dist}. Beginning with~\eqref{eq:nonchi} and averaging over $Z_0$, we can write the density function as follows
\begin{align}
f_{Z_t}\left(z_t\right)=&\int_0^\infty f_{Z_t|Z_0}\left(z_t|z_0\right)f_{Z_0}\left(z_0\right)\mathrm{d}z_0.
\end{align}
Then, letting $Z_0\sim\mathsf{Exp}\left(\theta\right)$, i.e., the initial condition of the process is drawn according to the  limiting distribution, we obtain
\begin{align}\nonumber
f_{Z_t}\left(z_t\right)=&\int_0^\infty ce^{-(z_tc+z_0u)}\mathrm{I}_0\left(2\sqrt{z_tuz_0c}\right)\frac{1}{\theta}e^{-z_0/\theta}\mathrm{d}z_0\\\nonumber
\labelrel={a}&\ \frac{ce^{-z_tc}}{u\theta+1}\sum^{\infty}_{j=0}\left(\frac{z_tuc\theta}{u\theta+1}\right)^j\frac{1}{j!}\\\nonumber
=&\frac{c}{u\theta+1}e^{-\frac{z_tc}{u\theta+1}}\\
=&\frac{1}{\theta}e^{-z_t/\theta}, 
\end{align}
where~\eqref{a} follows from the series expansion of $\mathrm{I}_0\left(2\sqrt{z_tuz_0c}\right)$.
%%%%%%%%%%%%%%%%%%%%%%%%%%%%%%%%%%%%%%%%
%%%%%%%%%%%%%%%%%%%%%%%%%%%%%%%%%%%%%%%%
\section{Proof of the autocovariance function of the squared distance process}
\label{Appen-Dist-acf}
In the following, we calculate the stationary autocovariance function of the process $\{Z_t, t\geq0\}$. In this case, the autocorrelation function depends only on the time shift $s-t=\Delta t$, i.e., 
\begin{align}\label{eq:autocovfor}
K_{Z}\left(\Delta t\right)&=\mathbb{E}\left\{Z_{s}Z_{t}\right\}-\mu_Z^2,
\end{align}
where $\mu_Z=\theta$ is the stationary mean. Substituting~\eqref{eq:SDEsol} in~\eqref{eq:autocovfor}, we obtain
\begin{equation}
K_{Z}\left(\Delta t\right)=\sigma^2e^{-k(t+s)}\mathbb{E}\left\{\left(\int_0^t\sqrt{Z_u}e^{ku}\mathrm{d}W_u\right)^2\right\}.
\end{equation}
Using the It\^{o} isometry $\left(\right.$Lemma $3.1.5$ in~\cite{oksendal2003stochastic}$\left.\right)$, the integral simplifies to
\begin{equation}
\mathbb{E}\left\{\left(\int_0^t\sqrt{Z_u}e^{ku}\mathrm{d}W_u\right)^2\right\}=\mathbb{E}\left\{\int_0^t Z_ue^{2ku}\mathrm{d}u\right\}.
\end{equation}
As $t,s \rightarrow\infty$, the stationary autocovariance function simplifies to
\begin{align}\nonumber
K_{z}\left(\Delta t\right)=\theta^2e^{-k\Delta t},
\end{align}
which has been shown in Proposition~\ref{prop:autocov}.
%%%%%%%%%%%%%%%%%%%%%%%%%%%%%%%%%%%%%%%%
%%%%%%%%%%%%%%%%%%%%%%%%%%%%%%%%%%%%%%%%
\section{Proof of the stochastic differential equation of the link SNR}
\label{Appen-SNR-SDE}
The squared distance process $Z_t$ is an It\^{o} process with stochastic differential equation given by~\eqref{eq:sde_z}. Now, let $N_t = g (t,Z_t)=Z_t^{-\eta/2}$. From It\^{o}'s formula $\left(\right.$Theorem $4.1.2$ in~\cite{oksendal2003stochastic}$\left.\right)$, $N_t$ is again an It\^{o} process with SDE
\begin{align}\nonumber
\mathrm{d}N_t=&\frac{\delta g}{\delta t}\left(t,Z_t\right)\mathrm{d}t+\frac{\delta g}{\delta z}\left(t,Z_t\right)\mathrm{d}Z_t\\
&+\frac{1}{2}\frac{\delta^2g}{\delta z^2}\left(t,Z_t\right)\left(\mathrm{d}Z_t\right)^2, 
\end{align}
where 
\begin{align}\nonumber
\frac{\delta g}{\delta t}\left(t,Z_t\right)&=0,\\\nonumber
\frac{\delta g}{\delta z}\left(t,Z_t\right)&=-\frac{\eta}{2}Z_t^{-\left(\eta/2+1\right)},\\\nonumber
\frac{\delta^2g}{\delta z^2}\left(t,Z_t\right)&=\frac{\eta}{2}\left(\frac{\eta}{2}+1\right)Z_t^{-\left(\eta/2+2\right)},\\\nonumber
\left(\mathrm{d}Z_t\right)^2&=\sigma^2Z_t\mathrm{d}t.
\end{align}
%%%%%%%%%%%%%%%%%%%%%%%%%%%%%%%%%%%%%%%%
%%%%%%%%%%%%%%%%%%%%%%%%%%%%%%%%%%%%%%%%
\section{Proof of the bivariate distribution of the link SNR}
\label{Appen-joint-cdf-snr}
Here, we provide a brief outline for the derivation of the bivariate cumulative distribution function in~\eqref{eq:joint_cdf_snr}.
\begin{align}\nonumber
F_{N_s,N_t}\left(\rho_s,\rho_t\right)&=\mathbb{P}\left(N_{s}\leq \rho_{s},N_{t}\leq \rho_{t}\right)\\\nonumber
&=\mathbb{P}\left(Z_{s}\leq \rho_s^{-2/\eta},Z_{t}\leq\rho_s^{-2/\eta}\right)\\\label{eq:int_joint-dist}
&=\int_{z_t=0}^{\rho_t^{-2/\eta}}\int_{z_s=0}^{\rho_s^{-2/\eta}}f_{Z_s|Z_t}f_{Z_t}\mathrm{d}z_s\mathrm{d}z_t.
\end{align}
Substituting~\eqref{eq:nonchi} and~\eqref{eq:sta_distrib_dist} in~\eqref{eq:int_joint-dist}, we obtain
\begin{align}\nonumber
F_{N_s,N_t}\left(\rho_s,\rho_t\right)\labelrel={a}&\frac{c}{\theta}\int_{z_t=0}^{\rho_t^{-2/\eta}}e^{-cz_t}\sum^{\infty}_{j=0}\frac{u^jc^jz_t^j}{j!\Gamma\left(j+1\right)}\\\nonumber
&\times\int_{z_s=0}^{\rho_s^{-2/\eta}}e^{-cz_s}z_s^j\mathrm{d}z_s\mathrm{d}z_t\\\nonumber
=&\frac{1}{\theta}\sum^{\infty}_{j=0}\frac{u^j}{c^{j+1}j!\Gamma\left(j+1\right)}\\\nonumber
&\times\int_{z_s=0}^{\rho_s^{-2/\eta}}e^{-cz_s}z_s^j\mathrm{d}z_s\int_{z_t=0}^{\rho_t^{-2/\eta}}e^{-cz_t}z_t^j\mathrm{d}z_t,
\end{align}
where~\eqref{a} follows from the series expansion of $\mathrm{I}_0\left(2\sqrt{z_suz_tc}\right)$.
%%%%%%%%%%%%%%%%%%%%%%%%%%%%%%%%%%%%%%%%
%%%%%%%%%%%%%%%%%%%%%%%%%%%%%%%%%%%%%%%%
\section{Proof of the density function of the link SNR}
\label{Appen-MT}
Here, we provide a brief outline for the derivation of the pdf expression in~\eqref{eq:pdf-snr-eta}. The conditional density function of $ N_t$, given $\Upsilon_t=\upsilon$, is equal to
\begin{equation}\label{eq:cond_distrib}
f_{ N_t|\Upsilon_t}\left( \rho|\upsilon\right)=\upsilon e^{- \rho \upsilon}.
\end{equation}
The unconditional distribution of the SNR results from marginalizing~\eqref{eq:cond_distrib} over the random variable $\Upsilon_t$, i.e.,  
\begin{equation}\label{eq:pdf_snr1}
f_{ N_t}\left( \rho\right)=\int_{0}^\infty f_{ N_t|\Upsilon_t}\left( \rho|\upsilon\right)f_{\Upsilon_t}(\upsilon)\ \mathrm{d}\upsilon.
\end{equation}
By the change-of-variables formula, the density function of $\Upsilon_t$ is given by
\begin{equation}\label{eq:pdf-ups}
f_{\Upsilon}\left(\upsilon\right)=\frac{2}{\eta\theta}\upsilon^{2/\eta - 1}e^{-\frac{\upsilon^{2/\eta}}{\theta}}.
\end{equation}
It follows that the stationary probability density of the process $\{N_t, t\geq0\}$ is equal to
\begin{align}
f_{N_t}\left(\rho\right)=\frac{2}{\eta\theta}\int_{0}^\infty\upsilon^{\frac{2}{\eta}}e^{-\rho \upsilon}e^{-\frac{\upsilon^{2/\eta}}{\theta}}\mathrm{d}\upsilon.
\end{align}
Given $\eta = p/q$ with p and q integers, we define
\begin{equation}\label{eq:Mellin_I}
I(x)=\int_{0}^\infty t^{2q/p}e^{-x t}e^{-\frac{t^{2q/p}}{\theta}}\mathrm{d}t. 
\end{equation}
We use the Mellin-transform (MT) method for the exact calculation of the integral $I(x)$. We can get the Mellin transform as
\begin{multline}
\mathcal{M}\left\{I(x);s\right\}=\int_0^{\infty} t^{2q/p}e^{-\frac{t^{2q/p}}{\theta}}\\
\times\left(\int_0^{\infty}x^{s-1}e^{-x t}\mathrm{d}x\right)\mathrm{d}t.
\end{multline}
Next, we make the substitution $u=xt$ in the internal integral and obtain
\begin{equation}
\mathcal{M}\left\{I(x);s\right\}=\frac{p}{2q}\theta^{\frac{2q+p-ps}{2q}}\Gamma(s) \Gamma\left(\frac{2q+p-ps}{2q}\right),
\end{equation}
for $0<\mathsf{Re}\left(s\right)<\frac{2q+p}{p}$. 

Then the inverse transform can be written as
\begin{multline}
I(x)=\frac{\theta^{\frac{2q+p}{2q}}}{2\pi i}\frac{p}{2q}\int_{\delta'-i\infty}^{\delta'+i\infty}\left(x\theta^{\frac{p}{2q}}\right)^{-s}\\
\times\Gamma(s)\Gamma\left(1+\frac{p}{2q}-\frac{ps}{2q}\right)\ \mathrm{d}s.
\end{multline}
Given $p,q \in \mathbb{N}$, we make the substitution $s=2qu$ and write
\begin{multline}\nonumber
I(x)= \frac{1}{2\pi i}\frac{p^{\frac{3}{2}+\frac{p}{2q}}}{\sqrt{2q}}\left(2\pi\right)^{1-q-\frac{p}{2}}\theta^{1+\frac{p}{2q}}\\\nonumber
\times\int_{\delta'-i\infty}^{\delta'+i\infty}\left(\left(\frac{x}{2q}\right)^{2q}\left(p\theta\right)^{p}\right)^{-u}\prod_{n=0}^{2q-1}\Gamma\left(u+\frac{n}{2q}\right)\\
\times\prod_{n=0}^{p-1}\Gamma\left(\frac{n+1}{p}+\frac{1}{2q}-u\right)\ \mathrm{d}u\\\nonumber
=\frac{p^{\frac{3}{2}+\frac{p}{2q}}}{\sqrt{2q}}\left(2\pi\right)^{1-q-\frac{p}{2}}\theta^{1+\frac{p}{2q}}\\
\times\MeijerG*{2q}{p}{p}{2q}{\frac{2pq-2q-p}{2pq}, \dots, \frac{-p}{2pq}}{0, \dots, \frac{2q-1}{2q}}{\left(\frac{x}{2q}\right)^{2q}\left(\frac{p}{c_0}\right)^{p}},
\end{multline}
where $G\left(\cdot\right)$ denotes Meijer's G furcation, $0<\delta'<\frac{2q+p}{p}$, and $(a)$ holds from the multiplication theorem.
%%%%%%%%%%%%%%%%%%%%%%%%%%%%%%%%%%%%%%%%
%%%%%%%%%%%%%%%%%%%%%%%%%%%%%%%%%%%%%%%%
\IEEEtriggeratref{20}
\bibliographystyle{ieeetr}
%\IEEEtriggeratref{24}
\bibliography{IEEEabrv,biblio}
\end{document}